\documentclass[12pt]{article}

\usepackage{latexsym}
\usepackage{amsopn,amscd}
\usepackage{amsfonts}
\usepackage{amssymb}
\usepackage{amsthm}
\usepackage{amsmath}
\usepackage{epsfig,graphics}
\usepackage{a4wide}
\usepackage[square,authoryear]{natbib}

\numberwithin{equation}{section}

\begin{document}

\title{\bf Origins and Breadth of the Theory of Higher Homotopies}

\author{
J.~Huebschmann$^1$
\\[0.3cm]
$^1$ USTL, UFR de Math\'ematiques\\
CNRS-UMR 8524
\\
59655 Villeneuve d'Ascq C\'edex, France\\
Johannes.Huebschmann@math.univ-lille1.fr
 }

\maketitle

\begin{abstract}

\footnotesize{

\noindent Higher homotopies are nowadays playing a prominent role
in mathematics as well as in certain branches of theoretical
physics. The purpose of the talk is to recall some of the
connections between the past and the present developments. Higher
homotopies were isolated within algebraic topology at least as far back as
the 1940's. Prompted by the failure of the Alexander-Whitney
multiplication of cocycles to be commutative, Steenrod  developed
certain operations which measure this failure in a coherent
manner. Dold and Lashof extended  Milnor's classifying space
construction to associative $H$-spaces, and a careful examination
of this extension led Stasheff to the discovery of $A_n$-spaces
and $A_{\infty}$-spaces as notions which control the failure of
associativity in a coherent way so that the classifying space
construction can still be pushed through.

Algebraic versions of higher homotopies have, as we all know, led
Kontsevich eventually to the proof of the formality conjecture.
Homological perturbation theory (HPT), in a simple form first
isolated by Eilenberg and Mac Lane in the early 1950's, has
nowadays become a standard tool to handle algebraic incarnations
of higher homotopies. A basic observation is that higher homotopy
structures behave much better relative to homotopy than strict
structures, and HPT enables one to exploit this observation in
various concrete situations which, in particular,
leads to the effective calculation of various invariants
which are otherwise intractable.

Higher homotopies abound but they are rarely recognized explicitly
and their significance is hardly understood; at times, their
appearance might at first glance even come as a surprise, for
example in the Kodaira-Spencer approach to deformations of complex
manifolds or in the theory of foliations.
 }
\end{abstract}

\newpage

\tableofcontents

\section{Introduction}

It gives me great pleasure to join in this celebration of Murray
Gerstenhaber's 80'th and Jim Stasheff's 70'th birthday. I had the
good fortune to get into contact with Jim some 25 years ago. In
1981/82 I spent six months at the Swiss Federal Institute of
Technology (Z\"urich) as a Research Scholar. At the time, I
received a letter from Jim asking for details concerning my
application of {\em twisting cochains\/} to the calculation of
certain group cohomology groups. What had happened? At Z\"urich, I
had lectured on this topic, and Peter Hilton was among the
audience. This was before the advent of the internet; not even
e-mail was available, and people would still write ordinary snail
mail letters. Peter Hilton travelled a lot and in this way
transmitted information; in particular, he had told Jim about my
attempts to do these calculations by means of twisting cochains.
By the way, since Peter Hilton was moving around some much, once
someone tried to get hold of him, could not manage to do so, and
asked a colleague for advice. The answer was: {\em Stay where you
are, and Peter will certainly pass by.\/}

At that time I knew very little about higher homotopies, but over
the years I have, like many of us, learned much from Jim's
insight, his habit of bringing his readers, students, and
coworkers out from \lq\lq behind the cloud of unknowing\rq\rq, to
quote some of Jim's own prose in his thesis. All of us have
benefited from Jim's generosity with ideas.

I cannot reminisce indefinitely, yet I would like to make two more
remarks, one related with language and in particular with language
skills: For example, I vividly remember, in the fall of 1987,
there was a crash at Wall Street. I inquired via e-mail---which
was then available---, whether this crash created a problem, for
Jim or more generally for academic life. His answer sounded
somewhat like \lq\lq Not a problem, but quite a tizzy here\rq\rq.
So I had to look up the meaning of \lq\lq tizzy\rq\rq\ in the
dictionary. This is just one instance of how I and presumably many
others profitted from Jim's language skills. Sometimes Jim answers
an e-mail message of mine in Yiddish--apparently his grandfather
spoke Yiddish to his father. There is no standard Yiddish spelling and,
when I receive such a message, to uncover it, I must read it aloud
myself to understand the meaning, for example \lq\lq OY VEH\rq\rq\
which, in standard German spelling would be \lq\lq Oh Weh\rq\rq.

I feel honoured  by the privilege to have been invited to deliver
this tribute talk. I would like to make a few remarks related to
 Murray Gerstenhaber. I have met Murray some 20
years ago when I spent some time at the Institute in Princeton.
From my recollections, Murray was then a member of the alumni board
of the Institute and was always very busy. We got into
real scientific and personal contact only later. In particular, I
was involved in reviewing some of the Gerstenhaber-Schack results,
and I will never forget that I learnt from Murray about Wigner's
approach to the idea of contraction. Also from time to time,
beyond talking about mathematics, we talked about history. For
example, Ruth Gerstenhaber once observed how people would gather for tea in the
Fuld Hall common room in the afternoon as usual around the table,
and no-one would say a word but, one after another, would
eventually leave the room murmuring \lq\lq There is no
counter-example.\rq\rq\ The perception of a mathematician through
a non-mathematician is sometimes revealing.

Before I go into the mathematical details of my talk, let us wish
many more years to Jim and Murray and their wives.

Let me now turn to my talk. There would be much more to say than
what I can explain in the remaining time. I shall touch on various
topics and make a number of deliberate choices and I will make the
attempt to explain some pieces of mathematics. However, my
exposition will be far from being complete or systematic and will
unavoidably be biased. For example there are higher homotopies
traditions in Russia and in Japan related with Lie loops, Lie
triple systems and the like which I cannot even mention, cf. e.~g.
\cite{kikkaone} and \cite{sabimikh}. There is a good account of
Jim Stasheff's contributions up to his 60'th birthday, published
at the occasion of this event \cite{mccleary}. This was just
before the advent of Kontsevich's proof of the formality
conjecture. I will try to complement this account and can thereby,
perhaps, manage to avoid too many repetitions. Also I will try to
do justice to a number of less well known developments.

\section{The formality conjecture}

Let me run right into modern times and right into our topic:
Algebraic versions of higher homotopies have, as we all know, led
Kontsevich eventually to the proof of the formality conjecture
\cite{kontssev}: Let $M$ be a smooth manifold, let
$A=C^{\infty}(M)$ and $L=\mathrm{Vect}(M)$, and consider  the
exterior $A$-algebra $\Lambda_AL$ on $L$. Let $\mathrm{Hoch}(A)$
denote the Hochschild complex of $A$, suitably defined, e.~g. in
the Fr\'echet sense. Given the vector fields $X_1,\dots,X_n$ on
$M$, let $\Phi_{X_1,\ldots, X_n}$ be the Hochschild cochain given
by
\[
\Phi_{X_1, \ldots, X_n}(a_1,\ldots, a_n) = \frac
1{n!}\sum\mathrm{sign}(\sigma)\prod_{j=1}^n
X_{\sigma(j)}(a_j),\quad a_1,\dots a_n \in A.
\]
By a version of a classical result of
Hochschild-Kostant-Rosenberg, the obvious map
\begin{equation}
\Lambda_AL \longrightarrow \mathrm{Hoch}(A),\quad X_1\wedge
\ldots\wedge X_n \mapsto \Phi_{X_1, \ldots, X_n}, \label{sh}
\end{equation} is an isomorphism on cohomology. That is to
say, the Hochschild cohomology of $A=C^{\infty}(M)$ amounts to the
graded algebra $\Lambda_AL$ of multi vector fields on $M$.

The standard Schouten-Nijenhuis bracket $[\,\cdot\, , \,\cdot\,]$
turns the suspension $s(\Lambda_AL)$ of $\Lambda_AL$--this is
$\Lambda_AL$, regraded up by 1, into an ordinary graded Lie
algebra. Here the grading convention is the standard one in
algebraic topology to the effect that, in particular, a
differential {\em lowers\/} degree by 1. Likewise, the familiar
Gerstenhaber bracket on $\mathrm{Hoch}(A)$ turns the suspension
$s(\mathrm{Hoch}(A))$ of $\mathrm{Hoch}(A)$ into an ordinary
differential graded Lie algebra. However, the morphism \eqref{sh},
while certainly being compatible with the differentials,  is {\em
not\/} compatible with the Lie brackets.

For any differential graded Lie algebra $\mathfrak g$, the
familiar C(artan) C(hevalley) E(ilenberg)-construction $\mathrm
S'[\mathfrak g]$ furnishes a d(ifferential) g(raded) coalgebra. In
fact, given $\mathfrak g$, differential graded Lie algebra
structures on $\mathfrak g$ can be characterized in terms of dg
coalgebra structures on the symmetric coalgebra $\mathrm
S'[s(\mathfrak g)]$ on the suspension $s(\mathfrak g)$ of
$\mathfrak g$: They correspond precisely to the dg coalgebra
structures determined by a linear term, the differential, and a
quadratic term, the bracket. This allows for immediate
generalization: An sh-Lie algebra is a vector space $\mathfrak g$
together with a coalgebra differential on the symmetric coalgebra
$\mathrm S'[s(\mathfrak g)]$ on the suspension $s(\mathfrak g)$ of
$\mathfrak g$. The {\em formality conjecture\/}, as formulated and
established by Kontsevich \cite{kontssev}, says that \eqref{sh}
extends to a {\em Lie algebra twisting cochain\/}
\begin{equation}
\tau\colon \mathrm S'[s^{2}(\Lambda_AL )] \longrightarrow
s(\mathrm{Hoch}(A)). \label{sh2}
\end{equation}
Here $\tau$ being a {\em twisting cochain\/} means that $\tau$
satisfies the {\em deformation\/} or {\em Maurer-Cartan\/}
equation. Such a Lie algebra twisting cochain furnishes an sh-map
from the ordinary (differential) graded Lie algebra $s(\Lambda_AL
)$ to the ordinary differential graded Lie algebra
$s(\mathrm{Hoch}(A))$.

The twisting cochain $\tau$ has homogeneous constituents $\tau_j$,
$\tau_1$ being essentially the above morphism \eqref{sh}. The
higher terms $\tau_j$ ($j\geq 2$) are an instance of higher
homotopies, and $\tau$ is an instance of an sh-map, a term created
by Jim Stasheff, inspired by terminology introduced by Sugarawa
\cite{sugawtwo}, see Section \ref{twenty} below; here \lq\lq
sh\rq\rq\ stands for \lq\lq strongly homotopic\rq\rq. Thus,
without having the language and notation of {\em higher
homotopies\/} and that of {\em deformations\/} at his
disposal---remarkably, both Murray Gerstenhaber and Jim Stasheff
are behind the scene at this point and both from 1963---,
Kontsevich would not even have been able to {\em phrase\/} the
formality conjecture. This confirms a variant of an observation
which, with a grain of salt, reads thus: {\em Mathematics consists
in continuous and discreet development of language and
notation.\/}

A key observation, advocated by Jim Stasheff from early on, is
this: Even though we start with strict objects, an sh-map between
them may lead to new insight, not necessarily available from
ordinary strict maps. This kind of observation has been
successfully exploited in rational homotopy theory for decades.
Kontsevich noticed its significance in an area at first
independent of rational homotopy and, furthermore, managed to
exhibit a particular sh-map which establishes the formality
conjecture.

R. Thom had raised the issue of existence of a graded commutative
differential graded algebra of cochains on a space \cite{thomtwo}.
This prompted the development of rational homotopy, starting
notably with  D. Quillen \cite{quilltwo} and D. Sullivan
\cite{sullitwo}. A space whose rational (or real) cochain algebra
is sh-equivalent to its cohomology algebra is said to be {\em
formal\/}, the term {\em formal\/} referring to the fact that the
rational homotopy type is then a {\em formal\/} consequence of the
structure of the cohomology ring. The term {\em formality
conjecture\/} derives from this tradition.

The statement of the formality conjecture implies, as we know,
that every Poisson bracket on a smooth manifold admits a
deformation quantization.

\section{Early History}

One of the origins of homotopy is Gau\ss' analytic expression for
the linking number of two closed curves (1833). One of the origins
of {\em higher homotopies\/} is the idea of a classifying space;
this idea goes again back to Gau\ss\  (1828). Another origin of
higher homotopies is the usage of resolutions. It is a common
belief, perhaps, that resolutions go back at least to Hilbert's
exploration of syzygies \cite{hilbert}. Hilbert studied syzygies
in order to show that the generating function for the number of
invariants of each degree is a rational function. He also showed
that, for a homogeneous ideal $I$ of a polynomial ring $S$, the
\lq\lq number of independent linear conditions for a form
of degree $d$ in $S$ to lie in $I$\rq\rq\ is a polynomial function
of $d$. However, this is not the entire story. The problem of
counting the number of conditions had already been considered for
some time; it arose both in projective geometry and in invariant
theory. A general statement of the problem, with a clear
understanding of the role of syzygies--but without the word,
introduced a few years later by Sylvester (1814--1897)
\cite{sylvester}--is given by Cayley (1821--1895) \cite{cayley}.
In fact, in a sense, Cayley somewhat develops what is nowadays
referred to as the {\em Koszul} resolution \cite{koszuone} more
than 100 years before Koszul. The terminology {\em homotopy\/} was
apparently created by H. Poincar\'e (1895). Poincar\'e also
introduced the familiar loop composition.  Thus we see that, in
the historical perspective, Jim Stasheff is in excellent company.

\section{Various 20'th century higher homotopies}
\label{twenty} Prompted by the failure of the Alexander-Whitney
multiplication of cocycles to be commutative, Steenrod  developed
the system of $\cup_i$-products \cite{Steenrod}. These induce the
squaring operations which, in turn, measure this failure of
commutativity in a coherent manner. The non-triviality of these
operations implies in particular that, over the integers, there is
no way to introduce a differential graded commutative algebra of
cochains on a space. The $\cup_i$-products entailed the
development of s(trongly)h(omotopy)c(ommutative) structures as
well as that of Steenrod operations.

An $A_{\infty}$-structure may be described as a system of higher
homotopies together with suitable coherence conditions. Massey
products \cite{Massey} may be seen as invariants of certain
$A_{\infty}$-structures. An elementary example arises from the
familiar Borromean rings, consisting of three circles which are
pairwise unlinked but all together are linked. The name \lq\lq
Borromean\rq\rq\ derives from their appearance in the coat of arms
of the house of the aristocratic Borromean family in Northern
Italy. If we regard these rings as situated in the 3-sphere, then
the cohomology ring of the complement is a trivial algebra, but
there is a Massey product of three variables detecting the
simultaneous linking of all three circles.

At the time Massey products were isolated, Jim Stasheff was a
graduate student at Princeton. His advisor J. Moore suggested he
look at the problem of determining when a cohomology class of a
based loop space $\Omega X$ was a {\em suspension\/} or a {\em
loop class\/}, i.~e. came from a cohomology class of $X$. In
pursuing this question, Stasheff was led to work of Sugawara
\cite{grouplike}, who had a recognition principle for
characterizing loop spaces up to homotopy type.

The ordinary loop multiplication on $\Omega X$ gives it the
structure of an H-space that is associative up to homotopy.
Moore's version of the loop space shows that there is a based loop
space which is homotopic to the familiar one for which the loop
multiplication is {\em strictly\/} associative. The conclusion is
that associativity is {\em not\/} a {\em homotopy invariant
property\/}; we owe Jim a complete understanding of the homotopy
invariance properties of associativity, and his solution furnishes
a clean recognition principle for loop spaces and, in fact, for an
entire hierarchy of spaces between loop spaces and H-spaces, the
loop spaces being spaces which admit a classifying space.

Specifically, Stasheff defined a nested sequence of homotopy
associativity conditions and called a space an $A_n$-space if it
satisfies the $n$'th condition. Every space is an $A_1$-space, an
H-space is an $A_2$-space, and every homotopy associative H-space
is $A_3$. An $A_{\infty}$-space has the homotopy type of a loop
space.

A.~Dold and R.~Lashof \cite{doldlash} generalized to associative
H-spaces Milnor's construction of a classifying space for a
topological group \cite{milnorut}. Jim Stasheff extended the
Dold-Lashof construction to $A_{\infty}$-spaces through his study
of homotopy associativity of higher order: an
$A_{\infty}$-structure precisely gives a classifying space. All
this was worked out in his thesis, published as \cite{stashone}.
Sugawara had introduced conditions for a group-like space, see the
definition in terms of the conditions 3.1--3.3 on p.~129 of
\cite{grouplike} to be imposed upon two maps related by what
Sugawara had called an iteration of the standard relations.
Altering the appropriate part of these conditions to suit the case
of associativity more precisely and naturally led Jim Stasheff,
apparently prompted by F. Adams, to isolating a now familiar
family of polyhedra, that of {\em associahedra}. We shall see
below that these polyhedra actually constitute an {\em operad\/}.
Moreover, following Sugawara \cite{sugawtwo}, Stasheff  defined
maps of $A_n$-spaces, referred to as $A_n$-maps, which are special
kinds of $H$-maps \cite{stashone} (Def. 4.4 p. 298); these maps
are homotopy multiplicative in a strong sense. Via Sugawara's
work, $A_n$-maps are related to the Dold and Lashof construction.
When the homotopies defining an $A_n$-map exist for all $n$, the
corresponding map is {\em strongly homotopy multiplicative\/} in
the sense of Sugawara \cite{sugawtwo} (p.~259). Thus the
sh-terminology we are so familiar with nowadays was born.

The algebraic analogue of an $A_n$-space in the category of
algebras is an $A_n$-algebra, the case $n={\infty}$ being included
here. The original and motivating example was provided by the
singular chains on the based loop space of a space. This notion,
and variants thereof, has found many applications. One such
variant, $L_{\infty}$-algebras, have already been mentioned. A key
observation here is that $A_{\infty}$-structures behave correctly
with respect to homotopy, which is not the case for strict
structures. What corresponds to the classifying space construction
in geometry is now the bar tilde construction. Inside the bar
tilde construction,  Massey products show up which determine the
differentials in the resulting bar construction spectral sequence.
Stasheff referred to these operations as {\em Yessam\/}
operations. History relates that once, at the end of a talk of
Jim's, S. Mac Lane asked the question: {\em Who was Yessam?}

Let me recall a warning, one of Jim's {\em favorite warnings\/} in
this context: When the differential of an $A_{\infty}$-algebra is
zero, the conditions force the algebra to be strictly associative
but there may still be non-trivial higher operations encapsulating
additional information, as the example of the Borromean rings
already shows where the non-triviality of the Massey product
reflects the triple linking.

Jim Stasheff continued to work in the realm of fibrations. There
is, for example, a notion of {\em topological parallel
transport\/} developed by him.  A recent joint article of J.
Stasheff and J. Wirth entitled {\em Homotopy transition
cocycles\/} \cite{staswirt} reworks and extends J. Wirth's thesis
written in 1965 under the supervision of J. Stasheff.

\section{Homological perturbations}

Homological perturbation theory (HPT) has nowadays become a
standard tool to construct and handle $A_{\infty}$-structures.
The term \lq\lq homological perturbation\rq\rq\ is apparently due to
J. Milgram. The
basic HPT-notion, that of {\em contraction}, was introduced in
Section 12 of \cite{EilenbergMacLane53}: A {\em contraction}
\begin{equation*}
   (X
     \begin{CD}
      \null @>\nabla>> \null\\[-3.2ex]
      \null @<<\pi< \null
     \end{CD}
    Y, h)
  \end{equation*}
consists of chain complexes $X$ and $Y$, chain maps $\nabla\colon
X \to Y $ and $\pi\colon Y \to X$, and a degree 1 morphism
$h\colon Y \to Y$ such that
\[
\pi\nabla = \mathrm{Id},\ \nabla \pi -\mathrm{Id} = dh+hd,\
h\nabla=0,\ \pi h = 0,\ hh=0.
\]
The notion of \lq\lq recursive structure of triangular
complexes\rq\rq\  in Section 5 \cite{Heller54} is also an example
of what was later identified as a perturbation. The \lq\lq
perturbation lemma\rq\rq\ is lurking behind the formulas in
Chapter II of Section 1 of \cite{Shi} and seems to have first been
made explicit by M. Barrat (unpublished). The first instance known
to us where it appeared in print is \cite{Brown}. Jim Stasheff
collaborated with various colleagues on questions related with
homological perturbation theory \cite{gugensta}, \cite{gulasta},
\cite{gulstatw} including myself \cite{huebstas}. An issue dealt
with in these papers, as well as in my joint paper
\cite{HuebschmannKadeishvili91} with T.~Kadeishvili, is that of
compatibility of the perturbation constructions with algebraic
structure. This issue actually shows up when one tries to
construct e.~g. models for differential graded algebras.

A homological algebra and higher homotopies tradition was
created as well by Be\-ri\-kash\-vili and his students in Georgia (at the
time part of the USSR). More precise comments about the historical
development until the mid eighties may be found in the article
\cite{HuebschmannKadeishvili91}, and some specific comments  about
the Georgian tradition in \cite{Huebschmann992}.

In the articles \cite{Huebschmann89}, \cite{Huebschmann89a},
\cite{Huebschmann89b}, I explored the compatibility of the
perturbation constructions with algebraic structure and developed
suitable algebraic HPT-constructions  to exploit
$A_{\infty}$-modules arising in group cohomology. In this vein, I
constructed suitable free resolutions from which I was able to do
explicit numerical calculations in group cohomology which until
today still cannot be done by other methods. In particular,
spectral sequences show up which do not collapse from $E_2$. These
results illustrate a typical phenomenon: Whenever a spectral
sequence arises from a certain mathematical structure, there is,
perhaps, a certain strong homotopy structure lurking behind, and the
spectral sequence is an invariant thereof. The
higher homotopy structure is then somewhat finer than the spectral
sequence itself.

\section{Quantum groups}

The issues of associativity and coassociativity, as clarified by
Jim Stasheff, play a major role in the theory of quantum groups
and variants thereof, e.~g. quasi-Hopf algebras. Suffice it to
mention here that Drinfel'd has introduced a notion of quasi-Hopf
algebra in which coassociativity of the diagonal is modified in a
way in which the pentagon condition plays a dominant role,
analogous to the hexagonal Yang-Baxter equation replacement for
cocommutativity. Now, given a quasi-Hopf algebra $A$, the
quasi-Hopf structure induces a multiplication $ \mathrm B \mathcal
C \times \mathrm B \mathcal C \longrightarrow \mathrm B \mathcal
C$ on the classifying space $\mathrm B \mathcal C$ of the category
$\mathcal C$ of $A$-modules, and the quasi-Hopficity says that
this multiplication is homotopy associative. More details and
suitable references  may be found in \cite{staseitt} and
\cite{stashele}.

\section{Operads}

The notion of $A_n$-spaces and the clarity they provide for the
recognition problem for topological groups became the basis for
the development of homotopy invariant algebraic structures. In
particular, the recognition problem for infinite loop spaces and
the simultaneous interest in coherence properties in categories
led to the idea of an operad \cite{maclasix}, \cite{maybotwo}.
With hindsight we recognize that a space is an $A_{\infty}$-space
if and only if it is an algebra over a suitably defined operad,
the non-symmetric operad $\mathcal K=\{K_n\}$ of {\em
associahedra\/}. In fact, this is the main result of Stasheff's
thesis, though not spelled out in this language:

{\em A connected space $Y$ of the homotopy type of a CW-complex
has the homotopy type of a loop space if and only if there exist
maps $K_n \times Y^n \to Y$ which fit together to make $Y$ an
algebra over the operad $\mathcal K$.} In fact, $Y$ then has the
homotopy type of the space $\Omega X$ where $X$ is constructed as
a quotient of $\coprod K_n \times Y^n$. This brings the
generalized classifying space construction to the fore.

Likewise a graded object is an $A_{\infty}$-algebra if and only if
it is an algebra over a suitably defined operad, and an
$L_{\infty}$-algebra can be characterized in the same manner as
well. In recent years many more new phenomena and structures and,
in particular, applications of operads have been found, in
particular in the theory of moduli spaces and in mathematical
physics.

\section{Deformation theory}

There is an obvious formal relationship between homological
perturbations and deformation theory but the relationship is
actually much more profound: In \cite{halpsttw}, Steve Halperin
and Jim Stasheff developed a procedure by means of which the
classification of rational homotopy equivalences inducing a fixed
cohomology algebra isomorphism can be achieved. Moreover, one can
explore the rational homotopy types with a fixed cohomology
algebra by studying perturbations of a free differential graded
commutative model by means of techniques from deformation theory.
This was initiated by M. Schlessinger and J. Stasheff
\cite{schlstas}. A related and independent development, phrased in
terms of what is called the functor $\mathcal D$, is due to N.
Berikashvili and his students at Georgia, notably T. Kadeishvili
and S. Saneblidze. Some details and references are given in
\cite{Huebschmann992}.
A third approach in which only the underlying graded vector space was
fixed is due to Y. Felix \cite{felix}.

More recently, prompted by a paper of Baranikov and Kontsevich
\cite{barakont}, Jim Stasheff and I developed an approach to
constructing solutions of the master equation by means of
techniques from HPT \cite{huebstas}. In that paper, we restricted
attention to contractions of a differential graded Lie algebra
onto its homology. More recently, I extended this approach to the
situation of a contraction of a differential graded Lie algebra
onto a general chain complex and thereby established the
perturbation lemma for differential graded Lie algebras
\cite{pertlie}. Further, I generalized the statement of the
perturbation lemma to arbitrary sh-Lie algebras \cite{pertlie2}.

\section {Strings}

Operads and sh-Lie algebras show up naturally in string and
conformal field theories, and Jim Stasheff contributed to this
area as well. Some details and more references may be found in
\cite{stasnint}.

\section{Cohomological physics}

One of Jim's long-term interests is physics. Due to his efforts it
is, perhaps, no longer a surprise that some structures of interest
in physics can be explored by means of tools going back to
topology, including graded Lie algebras and homological
perturbations. Jim contributed to {\em anomalies\/}
\cite{stashthr} and invested time and effort to unravel, for
example, the structure behind a field theory construction which
originally goes back to {\em Batalin, Fradkin, and Vilkovisky\/}.
The term \lq\lq cohomological physics\rq\rq\ was created by Jim.
See in particular \cite{stashnin} and \cite{stashsit} for details.

\section{Higher homotopies, homological perturbations, and the
working mathematician}

I have already explained how higher homotopies and homological
perturbations may be used to solve problems phrased in language
entirely different from that of higher homotopies and HPT. Higher
homotopies and HPT-constructions occur implicitly in a number of
other situations in ordinary mathematics where they are at first
not even visible. I can only mention some examples; these are
certainly not exhaustive.

\begin{itemize}

\item Kodaira-Nirenberg-Spencer: Deformations of complex
structures \cite{kodnispe};

\item Fr\"olicher spectral sequence of a complex manifold
\cite{Huebschmann98, Huebschmann00};

\item Toledo-Tong: Parametrix \cite{toltoone};

\item Fedosov: Deformation quantization \cite{fedosone};

\item Whitney, Gugenheim: Extension of geometric integration
to a contraction \cite{gugenhfi}, \cite{whitnsev}. Whitney's
geometric integration theory laid some of the ground work for
Sullivan's theory of rational differential forms quoted above. The
upshot of Gugenheim's contribution here is that the integration
map in de Rham theory is sh-multiplicative, the de Rham algebra
being an ordinary graded commutative algebra. This situation is
formally the same as that of the formality conjecture explained
above.

\item Huebschmann: Foliations \cite{omni}; in this paper, the
requisite higher homotopies are described in terms of a
generalized {\em Maurer-Cartan\/} algebra.

\item Huebschmann: Equivariant cohomology and Koszul duality
\cite{koszul}, \cite{koszultw}.

\item Operads; see e.~g. the conference proceedings which contain the article
\cite{stasnint}.
\end{itemize}

\end{document}